\documentclass[12pt]{amsart}
\usepackage{latexsym, amsmath,amssymb}
\usepackage{latexsym, amsmath,amssymb}
\usepackage{mathtools,slashed}
\usepackage{graphicx,subfigure,psfrag}
\usepackage{hyperref}
\usepackage{xcolor}
\hypersetup{
    colorlinks,
    linkcolor={red!50!black},
    citecolor={blue!50!black},
    urlcolor={blue!80!black}
}

\usepackage{graphicx}

 \numberwithin{equation}{section}

\def\XXint#1#2#3{{\setbox0=\hbox{$#1{#2#3}{%
\int}$ }
\vcenter{\hbox{$#2#3$ }}\kern-.6\wd0}}

\renewcommand{\epsilon}{\varepsilon}
\newtheorem{theorem}{Theorem}

\newtheorem{lemma}[theorem]{Lemma}
\newtheorem{corr}[theorem]{Corollary}

\newtheorem{proposition}[theorem]{Proposition}
\newtheorem{deff}[theorem]{Definition}

\newtheorem{conjecture}[theorem]{Conjecture}
\newcommand{\bth}{\begin{theorem}}
\newcommand{\ble}{\begin{lemma}}
\newcommand{\bcor}{\begin{corr}}

\newcommand{\bdeff}{\begin{deff}}

\newcommand{\bprop}{\begin{proposition}}
\newcommand{\ele}{\end{lemma}}
\newcommand{\ecor}{\end{corr}}
\newcommand{\edeff}{\end{deff}}

\numberwithin{theorem}{section}

\newcommand{\eprop}{\end{proposition}}

\renewcommand{\Pi}{\varPi}

\renewcommand{\epsilon}{\varepsilon}

\begin{document}

\title[Sign equidistribution of Legendre polynomials]
{Sign equidistribution of Legendre polynomials}

\author[A. D. Mart\'inez]{\'Angel D. Mart\'inez} \address{Albacete, Fields Institute, University of Toronto Mississauga, Canad\'a} \email{amartine@fields.toronto.edu}

\author[F. Torres de Lizaur]{Francisco Torres de Lizaur} \address{Sevilla, Fields Institute, University of Toronto, Canad\'a} \email{fj.torres@icmat.es}
\maketitle

\begin{abstract}
We prove {\em sign equidistribution} of Legendre polynomials:  the ratio between the lengths of the regions  in the interval $[-1, 1]$ where the Legendre polynomial assumes positive versus negative values, converges to one as the degree grows. The proof method also has application to the symmetry conjecture for a basis of eigenfunctions in the sphere.
\end{abstract}

\section{Introduction}

The importance of Legendre polynomials, from classical potential theory to modern computational methods, stems from the method of separation of variables in mathematical physics. They appear naturally in the spherical harmonic decompositions of functions in spherical coordinates. 

The zeros of these polynomials have been extensively studied in the past two centuries. They are known to be simple and belong to the interval $[-1,1]$. A classical result due to Bruns affirms that the roots $\theta_j$ of Legendre polynomials $P_n(\cos(\theta))$ equidistribute in $[0,\pi]$ as the degree $n$ grows. More concretely, if we denote the increasing sequence of zeroes  by
\[\theta_1<\theta_2<\cdots<\theta_n\]
the following inequalities hold
\[\frac{j-\frac{1}{2}}{n+\frac{1}{2}}\pi\leq\theta_{j}\leq\frac{j}{n+\frac{1}{2}}\pi\]
for $j=1,\ldots, n$. Markoff and Stieljes improved this to
\[\frac{j-\frac{1}{2}}{n}\pi\leq\theta_{j}\leq\frac{j}{n+1}\pi\]
for $j=1,\ldots,\lfloor n/2\rfloor$ that extends by symmetry considerations to inequalities for all the zeros (cf. \cite{S2}; or also the original articles \cite{M, S}). This was finally improved by Szeg\"o in 1936 who showed
\[\frac{j-\frac{1}{4}}{n+\frac{1}{2}}\pi\leq\theta_{j}\leq\frac{j}{n+1}\pi.\]

The main result of this paper explores yet another equidistribution property of Legendre polynomials, that we call sign equidistribution: we would say that a sequence of real polynomials $P_n$, or its zeros, sign equidistribute in an interval $I$ if the length of the set where the polynomial $P_n$ is positive equals the length of the set where the polynomial is negative in the limit $n\rightarrow \infty$. This notion was introduced in \cite{MT} in connection with the symmetry conjecture for the semiclassical limit of eigenfunctions on compact Riemannian manifolds.

\begin{theorem}[Sign equidistribution] \label{sign}
Let $\{\theta_{j}\}_{j=1}^{n}$ be the increasing sequence of zeros corresponding to the $n$th Legendre polynomial $P_n(\cos(\theta))$. For any closed interval $I\subseteq(0,\pi)$ containing an even number of roots we have
\[\bigg|\sum_{\theta_j\in I}(-1)^j\theta_j\bigg|=\frac{\text{ length }(I)}{2}+O(n^{-1})\]
where the constant is independent on $n$ but might depend on $I$.
\end{theorem}

Unfortunately, the bounds of Bruns-Szeg\"o are not enough to obtain this, and improving Szeg\"o's result seems a difficult task. We follow a different route; our method of proof provides a general result that is of independent interest. Indeed,

\begin{theorem}\label{general}
Let $\{\theta_{j}\}_{j=1}^{n}$ be the increasing sequence of zeros corresponding to the $n$th Legendre polynomial $P_n(\cos(\theta))$. Let $I\subseteq(0, \pi)$ be a fixed closed interval. For any function $f$ analytic in a neighbourhood of the interval $I$ the following holds
\[\sum_{\theta_j\in I}(-1)^j f(\theta_j)=\sum_{j} (-1)^{j} f \left(\frac{2\pi j-\pi/2}{2n+1}\right)+O(n^{-1})\|f\|_{L^{\infty}}\]
where the second sum extends over those $j$ such that $\theta_j\in I$ and the constant is independent on $n$ but might depend on $I$.
\end{theorem}

This result is intimately related to the so-called symmetry conjecture on the semiclassical limit of eigenfunctions. The symmetry conjecture asserts that on a given Riemannian manifold $(M,g)$ the area of positiveness of a Laplace-Beltrami eigenfunction tends to equal its area of negativeness as the eigenvalue grows. The conjecture has been disproved by the authors in \cite{MT}. The counterexamples are explicit but the proof is a computer assisted argument for the three dimensional flat torus. It is nevertheless easy to observe that the conjecture is true in the particular case of the two dimensional flat torus (loc. cit.). This might suggest its truth in the case of surfaces. 

In order to put the conjecture in context let us recall the following result contained in the seminal work of Donelly and Fefferman
\begin{theorem}[Corollary 7.10 in \cite{DF}]\label{donfef}
Let $(M,g)$ be a real analytic Riemannian manifold. There exists a constant $C$ such that, for any eigenfunction $\psi$ of the Laplace-Beltrami operator:
\[\frac{1}{C}\leq \frac{\operatorname{vol}(\{x\in M:\psi(x)>0\})}{\operatorname{vol}(\{x\in M:\psi(x)<0\})}\leq C.\]
\end{theorem}
We emphasize that the constant $C$ depends on the manifold, but not on the eigenvalue. This was improved to general smooth metrics on surfaces by Nadirashvili in \cite{N}. In the case of the sphere, it can be proved as a consequence of the Bruns-Szeg\"o inequalities. 

As a rather straightforward application of Theorem \ref{general} in the case of $f(z)=\cos(z)$ we will provide a partial result towards the symmetry conjecture in the two dimensional sphere:
\begin{conjecture}[Symmetry]\label{conj1}
Let $\{\psi_{\lambda}\}$ be a sequence of spherical harmonics. The limit
\[ \frac{\operatorname{vol}(\{x\in M:\psi_{\lambda}(x)>0\})}{\operatorname{vol}(\{x\in M:\psi_{\lambda}(x)<0\})}\rightarrow 1\]
holds as $\lambda$ grows to infinity.
\end{conjecture}
Before stating it let us introduce the set $\mathcal{B}$ of eigenfunctions on $\mathbb{S}^2$  that consists of the Legendre polynomials $P_n(\cos(\theta))$, and the eigenfunctions $P_n^m(\cos(\theta))\cos(m\varphi)$ and $P_n^m(\cos(\theta))\sin(m\varphi)$, where $P_n^m$ denotes the associated Legendre polynomials, $1\leq m\leq n$, $\varphi$ is the azymuthal angle variable and $\theta$ the polar angle variable. We emphasize that the linear combinations of these functions do not belong to $\mathcal{B}$ (otherwise it would simply contain all the spherical harmonics of degree $n$).

\begin{theorem}[Symmetry for a basis of eigenfunctions of $\mathbb{S}^2$]\label{basis}
For any sequence of eigenfunctions $\psi_n\in\mathcal{B}$ with increasing eigenvalue $n(n+1)$ :
\[\lim_{n\rightarrow\infty}\frac{\operatorname{vol}\{x\in\mathbb{S}^2:\psi_n(x)>0\}}{\operatorname{vol}\{x\in\mathbb{S}^2:\psi_n(x)<0\}}=1\]
\end{theorem}

We remark that in the case of tori of any dimension, the existence of a basis of eigenfunctions with the above property is trivial.

The paper is organized as follows.  In section \ref{proofgeneral} we present the proof of Theorem \ref{sign}, indicating the trivial changes that it requires to prove Theorem \ref{general}. Section \ref{proofbasis} is devoted to the proof of Theorem \ref{basis} as an application of these results.

\section{Proof of Theorem \ref{general}} \label{proofgeneral}

Before proceeding to the proof let us state a technical result we shall need later.

\begin{theorem}[Laplace formula]\label{asymptotic}
For any $\epsilon>0$, the asymptotic
\[P_n(\cos(\theta))=\sqrt{\frac{2}{n\pi\sin(\theta)}}\cos\left(\left(n+\frac{1}{2}\right)\theta-\frac{1}{4}\pi\right)+E(\theta)\,,\]
with $E(\theta)=O(n^{-\frac{3}{2}})$, holds uniformly for $\theta \in (\epsilon, \frac{\pi}{2}-\epsilon)$. The first derivative satisfies
\[\frac{\partial}{\partial\theta}P_n(\cos(\theta))=\frac{\partial}{\partial \theta}\sqrt{\frac{2}{n\pi\sin(\theta)}}\cos\left(\left(n+\frac{1}{2}\right)\theta-\frac{1}{4}\pi\right)+E'(\theta)\,,\]
with $E'(\theta)=O(n^{-\frac{1}{2}})$. The constants involved are independent on $n$ but do depend on the fixed $\epsilon>0$.
\end{theorem}
The first formula corresponds to the classical Laplace formula for which a number of proofs and refinements can be found in Szeg\"o's treatise \cite{S1}. The second part does not seem to have been noted explicitly in the literature but one can adapt the arguments in \cite{S1} to provide a proof. We provide details on the Appendix that complement the arguments within Szeg\"o's treatise.

The basic idea is to employ the argument principle of complex analysis for an specific choice of contour integration inside a strip containing $I\subset (\epsilon, \frac{\pi}{2}-\epsilon)$, which provides the identity
\[\frac{1}{2\pi i }\oint_{\Gamma}\frac{z}{P_n(\cos(z))} \frac{\partial P_n(\cos(z))}{\partial z} dz=\sum_{\theta_j \in I} (-1)^j\theta_j\]
where the contour has the form of a braid alternating winding number around consecutive zeros as in the figure. Notice that we can restrict our analysis to the $(0, \frac{\pi}{2})$, as the Legendre polynomials satisfy $P_{n}(-x)=(-1)^{n}P_{n}(x)$.  For the sake of clarity, let us focus now on the particular case $I:=(\epsilon, \frac{\pi}{2}-\epsilon)$.


\begin{figure}[t]
  \centering
   \psfrag{A}{$\frac{1}{n}$}
 \psfrag{B}{$\epsilon$}
 \psfrag{G}{$\Gamma$}
 \psfrag{P}{$\pi$}
 \psfrag{0}{$0$}

\includegraphics[scale=1,angle=0]{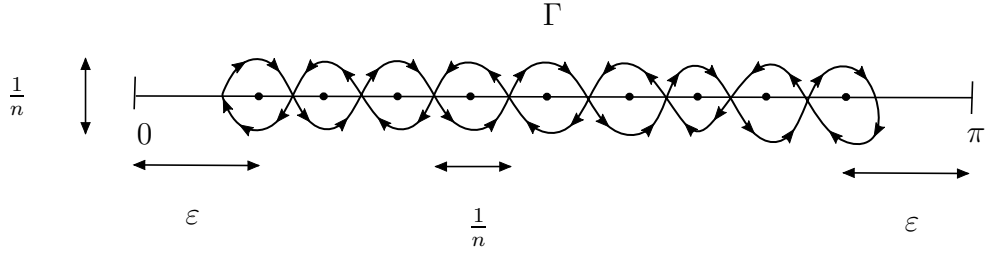}\hspace{0.1em}

\caption{The contour of integration $\Gamma$ braids around the roots at typical distance comparable to $n^{-1}$. The $j$th bounded region contains both $\theta_j$ and $\pi\frac{j-\frac{1}{4}}{n+\frac{1}{2}}$ (cf. Bruns-Sz\"ego's inequalities for the zeroes).}
    \label{fig:mesh1}

\end{figure}

We will denote
\[A(\theta):=\sqrt{\frac{2}{n\pi\sin(\theta)}}\cos\left(\bigg(n+\frac{1}{2}\bigg)\theta-\frac{1}{4}\pi\right)\]
to ease the notation. Theorem \ref{asymptotic} implies that the integral can be written as 

\[\frac{1}{2\pi i }\oint_{\Gamma}\left(\frac{z  A'(z)}{A(z)}+\frac{ z  E'(z)}{A(z)}- z (A'(z)+E'(z))\frac{E(z)}{A(z)(A(z)+E(z))}\right)dz\]
where we have used that
\[\frac{1}{A(z)+E(z)}-\frac{1}{A(z)}=\frac{-E(z)}{A(z)(A(z)+E(z))}.\]
We claim that the main term satisfies
\[\frac{1}{2\pi i }\oint_{\Gamma}\frac{z A'(z)}{A(z)}dz=\frac{\pi}{4}+O(n^{-1}).\]
To see this we only need to take into account the zeros of $\cos((n+1/2)\theta-\pi/4)$, given by
\[
\theta^{0}_{j}=\frac{j-\frac{1}{4}}{n+\frac{1}{2}}\pi,
\]
which clearly satisfy the sign equidistribution in $[0,\pi/2]$. To bound the remaining term, the idea is to show that it is equal to a gradient, plus some extra terms that go to zero as $n$ grows to infinity. More precisely, observe that, on the one hand
\[
z(A'+E')\frac{E}{A(A+E)}=z(A'+E')\frac{E}{A^2}-z(A'+E')\frac{E^2}{A^2(A+E)}
\]
and on the other hand,
\[
\frac{z E'}{A}-zA'\frac{E}{A^2}=\frac{\partial}{\partial z}\bigg(z \frac{E}{A}\bigg)-\frac{E}{A} \,.
\]
Putting both expressions together we see that the remaining term is equal to
\[
\frac{1}{2\pi i }\oint_{\Gamma} \frac{\partial}{\partial z}\bigg(z \frac{E}{A}\bigg) dz-\frac{1}{2\pi i }\oint_{\Gamma} \left(\frac{E}{A}+z(A'+E')\frac{E^2}{A^2(A+E)}-z E' \frac{E}{A^2}\right)dz.
\]
The first integral is clearly zero, since $z E/A$ is holomorphic in $\Gamma$. As for the other one, it can be bounded as 
\[
O\left(n^{-1}+\frac{n^{-3/2}}{\alpha^3}\|A'\|_{\infty}+\frac{n^{-1}}{\alpha^2}\right)\ell(\Gamma)=O(n^{-1}).
\]
Here $\ell(\Gamma)$ is the length of the contour $\Gamma$, and we have used the fact that $E=O(n^{-\frac{3}{2}})$, $E'=O(n^{-\frac{1}{2}})$ and, for the denominators, we claim $A=\Omega(n^{-\frac{1}{2}})$, i.e.  the fact that the contour can be chosen so that on it, $|A|\geq n^{-\frac{1}{2}} \alpha$ for some $\alpha>0$ depending on $\epsilon$ but independent of $n$. The bound $\|A'\|_{\infty}=O(n^{1/2})$ together with the above shows that
\[\sum_{\theta_j \in I}(-1)^j \theta_j=\frac{\pi}{4}+O(n^{-1})\]
as claimed.

To justify our claim let us consider the contour as in Figure \ref{fig:mesh1} that stays $\frac{1}{2n+1}$ away from the centered from the zeroes $\theta^{0}_j$ of $A(\theta)$. By Bruns-Sz\"ego inequality, the set of balls just described also contain the zeroes $\theta_j$ of $P_{n}$. On the other hand
\[
\alpha=\inf \left|\cos\left((n+\frac{1}{2})z-\frac{1}{4}\pi\right)\right|>0 \,
\]
where the infimum is taken on the complement to the union of balls, and it is independent of $n$.

The argument works, mutatis mutandis, for the zeroes contained in any other interval $I \subset [0, \pi]$, by simply adapting the contour $\Gamma$. Finally, Theorem \ref{general} can be proved verbatimly departing from
\[\frac{1}{2\pi i }\oint_{\Gamma}\frac{f(z) P'(\cos(z))\sin(z)}{P_n(\cos(z))}dz=\sum_{j}(-1)^jf(\theta_j)\]
instead. We leave details to the reader.

\section{Application to the symmetry conjecture}\label{proofbasis}

Let us observe first that any eigenfunction in the form of an associated Legendre polynomial already satisfies the conjecture, in fact (because of the symmetries of the $\cos(m \varphi)$ and $\sin(m \varphi)$ factors) the quotient is exactly one half for any degree! Furthermore, Legendre polynomials of odd degree  $k$ also satisfy the conjecture, since they verify $P_{k}(-z)=-P_{k}(z)$. Thus we can focus our attention to the Legendre polynomials of even degree $n$. 

The surface area of the part of the two dimensional sphere $C$ contained between two parallel planes $z=a$ and  $z=b$ is $2 \pi (b-a)$.

\begin{figure}[h]
    \centering
       \psfrag{a}{$z=a$}
 \psfrag{b}{$z=b$}

    \includegraphics[width=0.50\textwidth]{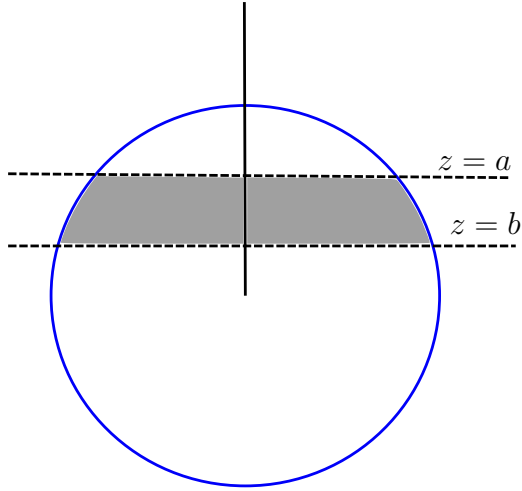}
    \caption{Schematic section of two dimensional spheres}
    \label{fig:mesh2}
\end{figure}

Thus, if $z_j$ are the roots of an even degree Legendre polynomial $P_n(z)$, either its area of positiveness or its area of negativeness is $2 \pi$ times the absolute value of the alternating sum
\[\sum_{j=1}^{n}(-1)^j z_j=\sum_{j=1}^{n}(-1)^j\cos(\theta_j).\]
Applying Theorem \ref{general} we have that, for any $\epsilon>0$, the above series is equal to
\[\sum_{\theta_j \in(\epsilon, \pi-\epsilon)} (-1)^j\cos\left(\frac{2\pi j-\pi/2}{2n+1}\right)+O_{\epsilon}(n^{-1/2})+O(\epsilon)\]
where the constants on the first might depend on $\epsilon$ but not the second (which is purely geometrical in nature, i.e. to compensate for the end points).

But on the other hand, it is easy to see that 
\[
I=\sum_{j=1}^{n} (-1)^j\cos\left(\frac{2\pi j-\pi/2}{2n+1}\right)=-1+ o(1) \,
\]
where $o(1)$ is a function that tends to zero as $n$ grows to infinity. To prove it we shall use the identity
\[
\cos\left(\frac{2\pi (j+1)-\pi/2}{2n+1}\right)-\cos\left(\frac{2\pi j-\pi/2}{2n+1}\right)=-\frac{2\pi}{2n+1} \sin\left(\frac{2\pi j-\pi/2}{2n+1}\right)+O(n^{-2})
\]
which comes from the Taylor approximation $\cos(x+\epsilon)=\cos(x)-\sin(x) \epsilon+O(\epsilon^2)$ valid uniformly in $[-\frac{\pi}{2},\frac{\pi}{2}]$. Using this we can rewrite the series as


\[\begin{aligned}
\sum_{j=1}^{n} (-1)^j\cos\left(\frac{2\pi j-\pi/2}{2n+1}\right)&=- \frac{1}{2} \sum_{k=0}^{\lfloor n/2 \rfloor} \frac{4\pi}{2n+1} \sin\left(\frac{2\pi (2k+1)-\pi/2}{2n+1}\right)+O(n^{-1})\\
&=-\frac{1}{2}\int_{0}^{\pi} \sin(\theta) d \theta+o(1)=-1+o(1) \,,
\end{aligned}\]
where in the last step we have recognized the sum as a Riemann integral approximation.

Summing up, we get
\[\sum_{j=1}^{n}(-1)^j z_j =-1 +o(1)+O_{\epsilon}(n^{-1/2})+O(\epsilon) \,.\]
 
Therefore, the difference 
\[\operatorname{vol}\{x\in\mathbb{S}^2:P_n(x)>0\}-\operatorname{vol}\{x\in\mathbb{S}^2:P_n(x)<0\}=O_{\epsilon}(n^{-1/2})+O(\epsilon).\]
One can divide by $\operatorname{vol}\{x\in\mathbb{S}^2:P_n(x)<0\}$, which is clearly bounded below away from zero, and then take the limit as $n$ grows, obtaining that
\[\lim_{n\rightarrow\infty}\frac{\operatorname{vol}\{x\in\mathbb{S}^2:P_n(x)>0\}}{\operatorname{vol}\{x\in\mathbb{S}^2:P_n(x)<0\}}=1+O(\epsilon) \,.\]

Since this is true for any $\epsilon>0$, this concludes the proof.


\section{Appendix}

As already mentioned the first asymptotic 
\[P_n(\cos(\theta))=\sqrt{\frac{2}{n\pi\sin(\theta)}}\cos\left((n+\frac{1}{2})\theta-\frac{1}{4}\pi\right)+E\]
with $E=O(n^{-3/2})$ is known as Laplace formula. The first part of the statement corresponds to Theorem 8.21.2 from \cite{S1}. We refer the reader to this reference for further details. To obtain the second part of Theorem \ref{asymptotic} one can follow the more general approach of Stieltjes (cf. Theorem 8.21.5 loc. cit.) which provides an error that can be explicitly written as
\[E(\theta)=\frac{2}{\pi}\operatorname{Im}\left(\frac{e^{i(n+1)\theta}e^{i(\pi/4-\theta/2)}}{(2\sin\theta)^{1/2}}\int_0^1t^n(1-t)^{-1/2}\frac{1}{\pi}\int_0^{\pi}\frac{z\sin^2(\varphi)}{1-z\sin^2(\varphi)}d\varphi dt\right)\]
where
\[z=(1-t)\frac{e^{i(\theta-\pi/2)}}{2\sin\theta}\]
cf. section 8.5 in \cite{S1}, specifically equation 8.5.1 considering $R_p(\theta)$ for $p=1$. Taking derivatives in the identity above and using the fact that $\theta\in(\epsilon,\pi-\epsilon)$, so that $\sin(\theta)$ is bounded away from zero, one concludes the proof as in the original (cf. bound 8.5.5 loc. cit.).

\end{document}